\documentclass[11pt]{amsart}
\usepackage[utf8]{inputenc}

\title{All graphs are majority 3-choosable}
\author{Jan Ouborny}
\author{Max Pitz}
\address{Universit\"at Hamburg, Department of Mathematics, Bundesstrasse 55 (Geomatikum), 20146 Hamburg, Germany}
\email{\{max.pitz, jan.ouborny\}@uni-hamburg.de}

\date{}

\usepackage{thmtools, thm-restate}%For duplicating theorem numbers.
\usepackage{amsmath,amssymb,amsthm}  
\usepackage{mathtools}
\usepackage{tikz}
\usetikzlibrary{decorations.pathreplacing, positioning}
\usepackage{comment}
\usepackage{enumitem}

\usepackage{todonotes}

\usepackage{xcolor} 	
\usepackage[unicode]{hyperref}
\hypersetup{
	colorlinks,
    linkcolor={red!60!black},
    citecolor={green!60!black},
    urlcolor={blue!60!black},
}
\usepackage[abbrev, msc-links]{amsrefs} 

\usepackage[utf8]{inputenc}
\usepackage[T1]{fontenc}
\usepackage{lmodern}
\usepackage[babel]{microtype}
\usepackage[english]{babel}

\usepackage{geometry}
\usepackage{enumitem}

\let\polishlcross=\l
\def\l{\ifmmode\ell\else\polishlcross\fi}

\let\emptyset=\varnothing

\let\theta=\vartheta
\let\rho=\varrho
\let\phi=\varphi

\def\NN{\mathbb N}

\def\cA{{\mathcal A}}

\def\cX{{\mathcal X}}
 
\newcommand{\parentheses}[1]{{\left( {#1} \right)}}

\newcommand{\p}{\parentheses}

\newcommand{\Set}[1]{{\left\lbrace {#1} \right\rbrace}}

\def\set#1:#2{\Set{{#1} \colon {#2}}}

\renewcommand{\subset}{\subseteq}

\newcommand{\saturated}{saturated}

\DeclareFontFamily{U}  {MnSymbolC}{}
\DeclareSymbolFont{MnSyC}         {U}  {MnSymbolC}{m}{n}
\DeclareFontShape{U}{MnSymbolC}{m}{n}{
    <-6>  MnSymbolC5
   <6-7>  MnSymbolC6
   <7-8>  MnSymbolC7
   <8-9>  MnSymbolC8
   <9-10> MnSymbolC9
  <10-12> MnSymbolC10
  <12->   MnSymbolC12}{}
\DeclareMathSymbol{\powerset}{\mathord}{MnSyC}{180}

\theoremstyle{plain}
\newtheorem{maintheorem}{Theorem}
\newtheorem{thm}{Theorem}[section]

\newtheorem{prop}[thm]{Proposition}

\newtheorem{lemma}[thm]{Lemma}

\newtheorem*{unfr}{Unfriendly Partition Conjecture}

\theoremstyle{definition}

\DeclareMathOperator{\nbly}{nbly}

\DeclareMathOperator{\dom}{dom}

\begin{document}
\begin{abstract}
Every graph is majority 3-choosable.  This generalises the result by Shelah-Milner that every graph has an unfriendly 3-partition, confirming a conjecture of Haslegrave from 2020. 
\end{abstract}

\maketitle

\section{Introduction}

A bipartition of the vertex set of a graph is \emph{unfriendly} if every vertex has at least as many neighbours in the other class as in its own. Every finite graph has an unfriendly partition: just take any bipartition that maximizes the number of edges between the partition classes. The corresponding problem for infinite graphs has become one of the best-known open problems in infinite graph theory. 

\begin{unfr}
    Every countable graph admits an unfriendly bipartition.
\end{unfr}

Shelah and Milner \cite{shelah1990graphs} constructed uncountable graphs that admit no unfriendly bipartitions. On the affirmative side, they proved that every infinite graphs has an unfriendly 3-partition, i.e.\ a partition of the vertex set into 3 partition classes such that every vertex has at least as many neighbours in the other two classes as in its own.

Unfriendly partitions may be rephrased in the language of colourings. A \emph{majority colouring} of a graph is an assignment of colours to vertices such that every vertex has at least as many neighbours of a colour different from its own colour than neighbours of the same colour. A graph is \emph{majority $\ell$-colourable} if it has a majority colouring using at most $\ell$ colours. %If a vertex has infinitely many neighbours, we interpret “at least half” to mean a set of the same cardinality as the whole neighbourhood. 
In this language, the Unfriendly Partition Conjecture asserts that every countable graph is majority 2-colourable, and Shelah and Milner's result means that every graph is majority 3-colourable.

This colouring perspective suggests the following list-version  of the unfriendly partition conjecture \cite{kreuzerMajority2017}: A graph $G$ is \emph{majority $\ell$-choosable} if for any system of lists $\set{L(v)}:{v \in V(G)}$ all of size $\ell$ there is a majority colouring of $G$ in which each vertex $v$ is coloured with a colour from its  list $L(v)$. 
Anholcer, Bosek and Grytczuk \cite{ANHOLCER2024103829} showed that every countable graph is majority $4$-choosable, and Haslegrave \cite{haslegrave2023} improved this by showing that every countable graph is majority $3$-choosable. Moreover, Haslegrave conjectured that the same result holds in fact for all graphs, also the uncountable ones \cite[Conjecture~4]{haslegrave2023}. This would be a natural strengthening of the Shelah-Milner result. The purpose of this note is to confirm this conjecture:

\begin{maintheorem}
\label{thm_intro}
    Every  graph is majority $3$-choosable.
\end{maintheorem}

Clearly, finite graphs are majority $2$-choosable by the same maximality argument as above. 
Thus, it remains to prove the case when $G$ is infinite.
Following the approach of Shelah-Milner, we establish the following technical strengthening of Theorem~\ref{thm_intro}, which lends itself better to an inductive proof.

Given a graph $G$, we write $N(v)$ for the neighbourhood of a vertex $v$, and $d(v) = |N(v)|$ for the \emph{degree} of $v$. Given a set of vertices $A \subset V(G)$, we abbreviate $d_A(v) = |N(v) \cap A|$.
We say a set of vertices $A$ in a graph $G$ is \emph{closed} if every vertex $v \in V \setminus A$ satisfies $d_A(v) < d(v)$. Note that the empty set is closed provided that $G$ contains no isolated vertices (which we may clearly assume).

A vertex $v$ of $G$ is \emph{happy} with a colouring of $G$ if at most half of the edges incident with $v$ are monochromatic.
\begin{maintheorem}
\label{thm_intro2}
     
    Let $G=(V,E)$ be an infinite graph and suppose each $v \in V$ has a list $L(v)$ of three colours.
    Suppose further that $V = A \sqcup B$ is a partition such that
    $A$ is closed and $B$ is infinite.
    Let $h \colon A \to C:= \bigcup_{v} L(v)$ be any list-colouring of $A$, and 
    fix $x \in B$ and $c_x \in L(x)$ be arbitrary. 
    Then there is a list-colouring $g \colon V \to  C$ extending $h$ such that every vertex of $B$ is happy with $g$ and such that $g(x) \neq c_x$.
\end{maintheorem}

 Then Theorem~\ref{thm_intro} follows from Theorem~\ref{thm_intro2} with $A=\emptyset$ and $B=V$. 
 We prove Theorem~\ref{thm_intro2} by induction on the cardinality of $B$. In the next section, we deal with the case when $B$ is countably infinite (following the strategy of Haslegrave), and in the last section with the case when $B$ is uncountable (following the approach of Shelah and Milner).

\section{The countable case}

The following lemma has been essentially observed by Anholcer, Bosek, and Grytczuk \cite{ANHOLCER2024103829}, and also by Haslegrave \cite[Lemma~1]{haslegrave2023}. We  verify that it also holds with the extra parameter $c_x \in L(x)$.

\begin{lemma} \label{lem_pre_processing}
    Let $V$ be a countably infinite set, and $\cX$ a countable (possibly finite) collection of infinite subsets of $V$. 
    Suppose that every $v \in V$ has a list $L(v)$ of $\ell + 1$ colours. 
    Further, let $x \in V$ and $c_x \in L(x)$ be arbitrary.
    
    Then there is a choice of $\ell$-element sublists $L'(v) \subsetneq L(v)$ for every $v \in V$, where $L'(x)= L(x) \setminus \Set{c_x}$, with the property that for every colour $c$ and set $X \in \cX$ there are infinitely many $v \in X$ such that $c \not \in L'(v)$.
\end{lemma}
\begin{proof}
    First, we note that $C\coloneqq \bigcup_{v \in V} L(v)$ is countable. 
    Without loss of generality, we assume that $V \in \cX$, so $\cX \neq \emptyset$.
   Then $\cX \times C \times \NN$ is countably infinite, and 
    we fix an enumeration of $\cX \times C \times \NN$ such that $(V, c_x, 1)$ appears as the first element.
    We also fix an enumeration of $V$ such that $x$ appears as the first element. 

    We choose for every triple $(X,c,n) \in \cX \times C \times \NN$ the minimal (regarding the enumeration of $V$) vertex $v \in X$ that has not been chosen before.
    Since $X$ is infinite, $v$ exists.
    If $c \in L(v)$, then we set $L'(v) \coloneqq L(v) \setminus \Set{c}$. 
    Otherwise, we define $L'(v)$ as an arbitrary sublist of $L(v)$ with $\ell$ elements.
    
    Since $x$ appears first in the enumeration of $V$ and $(V, c_x, 1)$ appears first in the enumeration of $\cX \times C \times \NN$, we have $L'(x)= L(x) \setminus {c_x}$.
    Finally, we note that for every pair $(X,c) \in \cX \times C$ every $v$ that was chosen in the construction step of $(X,c,n)$ with $n \in \NN$ has the property that $v \in X$ and $c \not \in L'(v)$.
    Since there are infinitely many distinct triples $(X,c,n)$ with $n \in \NN$, there are infinitely many such vertices. 
    Hence, the chosen sublists have the desired property.
\end{proof}

\begin{proof}[Proof of Theorem~\ref{thm_intro2} for countable $B$] 
    Let $b_1, b_2, \dots$ be an enumeration of the set $B$ with $b_1 = x$.
    Apply Lemma \ref{lem_pre_processing} on the countable graph $G[B]$ with $\cX = \Set{N(b_n) \cap B \colon  n \in \NN \text{ and } d_B(b_n) = \infty}$ and $\ell=2$ and $x,c_x$ to obtain a family of lists $(L'(b_n))_{n \in \NN}$ all of size $2$ with $c_x \not \in L'(x)$.

    For $n \in \NN$ put $B_n=\Set{b_1, \dots, b_n}$, and consider the subgraph $G_n$ of $G$ induced by $B_n \cup \bigcup_{i=1}^n (N(b_i) \cap A)$.
    Since $A \cap B = \emptyset$, $A$ is closed and $B$ is countable, each $b\in B$ has only finitely many neighbours in $A$.
    Therefore, each $G_n$ is finite.
    We claim that for each $n \in \NN$ there is a list-colouring $g_n$ of $G_n$ such that 
    \begin{enumerate}
        \item $g_n$ and $h$ agree on their common domain,
        \item $g_n(b) \in L'(b)$ for every $b \in B_n$,
        \item $g_n(b_1) \neq c_x$,
        \item every vertex of $B_n$ is happy with $g_n$.
    \end{enumerate}
    Let $g_n$ be a list-colouring of $G_n$ satisfying (1)--(3) and with a maximal number of cross edges induced by the colour classes of $g_n$. 
    Since $G_n$ is finite, there are only finitely many list-colourings that satisfy (1)--(3). Hence, $g_n$ exists, 
    and it is straightforward to check that the maximality of $g_n$ implies that it also satisfies (4).

    Now we continue with the compactness argument from \cite{haslegrave2023} to obtain our desired majority colouring $g \colon V \to C$. First, $g$ shall agree with $h$ on $A$. In order to extend $g$ to $b$, observe that
    infinitely many colourings of $(g_n)_{n \in \NN}$ agree on $b_1$; 
    we define $g(b_1)$ as this colour.
    Of these infinitely many colourings that agree on $b_1$, infinitely many agree on $b_2$; we define $g(b_2)$ as this colour. And so on. 
    We continue this construction to obtain $g \colon V \to C$ with $g(b) \in L'(b)$ for every $b \in B$.

    We claim that $g$ is our desired list-colouring. By construction, $g$ extends $h$, and satisfies $g(x) = g(b_1) \neq c_x$ by (3). It remains to show that every vertex $b \in B$ is happy with $g$.
    If $d_B(b) < \infty$, then put $k = \max \Set{\ell \colon b_\ell \in N(b) \cup \Set{b}}$.
    Then there exists $m \geq k$ such that $g_m(b_\ell)=g(b_\ell)$ for every $\ell \leq k$. 
    Since $b$ is happy with $g_m$, and $g_m$ agrees with $g$ on all neighbours of $b$ by (1), it follows that $b$ is happy with $g$.
    Otherwise,
    $d_B(b) = \infty$. Since $b$ has only finitely many neighbours in $A$ (as $A$ is closed), it has countable degree in $G$. So $b$ is happy provided it has infinitely many neighbours with different colour as its own. To this end, consider $X=N(b) \cap B \in \cX$ and $c=g(b)$. 
    By the definition of the sublists according to Lemma \ref{lem_pre_processing}, we get that $b$ has infinitely many neighbours $v \in X$ with $g(b) = c \not \in L'(v) \ni g(v)$. 
\end{proof}

\section{The uncountable case}

We begin with a lemma from Shelah-Milner, filling in some additional details in the proof.

\begin{lemma}[{\cite[Lemma 1]{shelah1990graphs}}]  \label{lemma_1}
Let $\cA = \set{A_i}:{i \in I}$  be an infinite family of sets such that $|A_i| \geq |I|$ for all $i \in I$. Then there are pairwise disjoint sets $B_i \subseteq A_i$ such that $|B_i| = |A_i|$.
\end{lemma}

\begin{proof}
Let $D = \set{|A_i|}:{ i \in I}$ list the cardinalities of the sets involved, and let
$$R = \set{\kappa \in D}:{ \kappa > \sup \set{\mu}:{\mu \in D, \; \mu < \kappa}}.$$
For $\kappa \in R$ let $I(\kappa) = \set{i \in I}:{|A_i| \geq \kappa}$. 
Recursively in $\kappa \in R$, we choose subsets
$$A_i (\kappa) \subset A_i \quad \text{for }i \in I(\kappa)$$
such that $|A_i(\kappa) | = \kappa$ and $A_i(\kappa) \cap A_j(\mu) = \emptyset$ whenever $(\kappa,i) \neq (\mu,j)$. 

Indeed, consider $\kappa \in R$ and suppose that we have already chosen all $A_i(\mu)$ for $\mu < \kappa$  accordingly. Let $X$ denote the union over all these previously chosen $A_i(\mu)$. Then $\kappa \in R$ implies $|X| < \kappa$. %
Pick a function $f \colon \kappa \to I(\kappa)$ such that $|f^{-1}(i)| = \kappa$ for all $i \in I(\kappa)$. 
For $\alpha < \kappa$, we now recursively choose $x_\alpha \in A_{f(\alpha)} \setminus \p{X \cup \set{x_\beta}:{ \beta < \alpha}}$, which is possible since $|A_{f(\alpha)}| \geq \kappa$ and $|X \cup \set{x_\beta}:{ \beta < \alpha}| < \kappa$. Then $A_i(\kappa) = \set{x_\alpha}:{\alpha \in f^{-1}(i) } \subseteq A_i$ is as desired.

At the end of the construction, we claim the sets
$$B_i = \bigcup \set{A_i(\kappa)}:{ i \in I(\kappa), \; \kappa \in R}$$
are as desired.
 Indeed, as all $A_i(\kappa)$ are pairwise disjoint, the $B_i$ are disjoint. Further, we have $B_i \subseteq A_i$. It remains to argue that $|B_i| = |A_i|$. If $\kappa=|A_i| \in R$, then $A_i(\kappa) \subseteq A_i$ witnesses that $|B_i| = \kappa$. Otherwise, if $\kappa=|A_i| \in D \setminus R$, then $\kappa = \sup \set{\mu}:{\mu \in D, \; \mu < \kappa}$ by definition of $R$. Since $\sup \set{\mu}:{\mu \in D, \; \mu < \kappa} = \sup \set{\mu}:{\mu \in R, \; \mu < \kappa}$,
 it follows that $\bigcup \set{A_i(\mu)}:{\mu < \kappa} \subset A_i$ has size $\kappa$ are desired.
\end{proof}

The following concept of ``saturation'' is implicit in the Shelah-Milner proof.
Let $A$ be a subset of $V(G)$ of a graph $G$.
We say $B^* \subseteq V(G)$ is $A$--\emph{\saturated} if for every $b \in B^*$ with $N(b) \setminus (A \cup B^*) \neq \emptyset$ we have 
$$|N(b) \setminus (A \cup B^*)| > |B^*| \quad \text{and} \quad |(N(b)\setminus A) \cap B^*| =|B^*|.$$

\begin{lemma}[Saturation lemma] \label{closure_lemma}
    Let $A$ be a prescribed set of vertices in a graph $G$. Then every infinite set of vertices $B$ in $G$ is included in an $A$--\saturated\ set of vertices  $B^*$ with $|B|=|B^*|$ and $B^* \cap A= B \cap A$.
\end{lemma}
\begin{proof}
    Let $\mu = |B|$ and put $B_0 \coloneqq B$. We construct $B_{n+1}$ from $B_n$ by adding for every $b \in B_n$ with $|N(b) \setminus A| \leq \mu$ the set $N(b) \setminus A$ to $B_n$; and for every $b \in B_n$ with $|N(b) \setminus A| > \mu$ we add an arbitrary $\mu$-sized subset of $N(b) \setminus A$ to $B_n.$

    We claim that $B^* \coloneqq \bigcup_{n \in \NN} B_n$ is as desired.
    First, we show by induction that $|B_n|=\mu$ and $B_n \setminus B \subseteq V \setminus A$ for every $n \in \NN$.
    Clearly, $|B_0|=\mu$ and $B_0 \setminus B_0 \subseteq V \setminus A$.
    Suppose that $|B_n|=\mu$ and $B_n \setminus B \subseteq V \setminus A$. Then there are at most $\mu$ many vertices in for which we added at most $\mu$ many vertices of $V \setminus A$ to obtain $B_{n+1}$.
    Hence, $|B_{n+1}|=\mu$ and $B_{n+1} \setminus B \subseteq V \setminus A$. In particular, we get that $B^* \setminus B \subseteq V \setminus A$, which is equivalent to  $B^* \cap A= B \cap A$.
    
    Next, let $b \in B^*$ with $N(b) \setminus (A \cup B^*) \neq \emptyset$. 
    Since $N(b) \setminus (A \cup B^*)$ is non-empty, we have $|N(b) \setminus A| > \mu$, 
    as otherwise we would have added $N(b) \setminus A$ to $B^*$. 
    With $|B^*|=\mu$, it follows that $|N(b) \setminus (A \cup B^*)| > \mu$.
    Finally, by the construction we have added a $\mu$-sized subset of $N(b) \setminus A$ to $B_{n+1}$. Hence, we have $|(N(b)\setminus A) \cap B^*|=\mu$. 
\end{proof}

Let $G=(V,E)$ be a graph. 
We define for every subset $A \subseteq V$ the set 
\[
    \nbly(A) = \Set{v \in V \colon d_A(v) = d(v)}.
\]
Then $A$ is closed if and only if $\nbly(A) \subseteq A$.
If $\set{A_i}:{i \in I}$ are closed subsets, then the intersection $A = \bigcap_{i \in I} A_i$ is closed as well. To see this, fix $v \in V \setminus A$. Then there is $i \in I$ with $v \notin A_i$. Since $A \subset A_i$, and $A_i$ is closed, it follows that
$d_A(v) \leq d_{A_i}(v) < d(v).$
By Zorn's lemma, this allows us to define the \emph{closure} $\overline{A}$ of $A$ as the unique, inclusion-wise minimal closed set containing $A$. Write $\partial A = \overline{A} \setminus A$ for the \emph{boundary} of $A$.

\begin{prop}\label{prop_enumerate_boundary} Let $A \subseteq V$. The following statements are true
    \begin{enumerate}[label=(\roman*)]
        \item % $|N(v)|=|N(v) \cap \overline{A}|$ 
        $d(v) = d_{\overline{A}}(v)$ for every $v \in \partial A$, and 
        \item there is a well-ordering $\Set{a_i \colon i < \lambda}$ of $\partial A$ such that
        \[
            |N(a_i)|=|N(a_i) \cap (A \cup \Set{a_j \colon j < i)}|
        \] for every $i < \lambda$.
    \end{enumerate}
\end{prop}
\begin{proof}
    First, we define a transfinite increasing sequence $\Set{A_\alpha \colon \alpha < |G|^+}$ and for every non-limit ordinal $\alpha < |G|^+$ we define a well-ordering $<_\alpha$ of $A_{\alpha +1 } \setminus A_\alpha$.
    Let $A_0 \coloneqq A$. We construct $A_{\alpha +1}$ from $A_\alpha$ by adding all vertices in $\nbly(A_\alpha)$. 
    We fix an arbitrary well-ordering of $A_{\alpha +1 } \setminus A_\alpha$ define $<_\alpha$ as that well-ordering.
    If $\gamma$ is a limit ordinal and $A_\alpha$ has been constructed for all $\alpha < \gamma$, then we set $A_\gamma \coloneqq \bigcup_{\alpha < \gamma} A_\alpha$.
    Then there exists a minimal ordinal $\sigma < |G|^+$ such that $A_{\sigma+1}=A_\sigma$. 
    
    We claim that $A_\sigma = \overline{A}$. By construction, we have $A \subseteq A_\sigma$ and $A_\sigma$ is closed.
    Suppose for a contradiction there is a closed $B$ with $A \subset B \subsetneq A_\sigma$. 
    Let $\beta \leq \sigma$ be minimal such that $A_\beta \not \subseteq B$
    and choose $v \in A_\beta \setminus B$.
    By construction there is $\gamma < \beta$ such that $v \in A_{\gamma+1} \setminus A_\gamma$ and $v \in \nbly(A_\gamma)$.
    Since $\gamma < \beta$, we have $A_\gamma \subseteq B$. 
    Since $B$ is closed, we have $\nbly(A_\gamma) \subseteq B$, but $v \not \in B$, a contradiction.

    Now towards (i). Let $v \in \partial A$. By the construction, there is $\beta < \sigma$ such that $v \in \nbly(A_\beta)$. Hence, we have $|N(v)|=|N(v) \cap A_\beta|$; and since $A_\beta \subseteq \overline{A}$, we have $|N(v)|=|N(v) \cap \overline{A}|$.
    
    Finally, we consider the well-ordering $\Set{a_i \colon i < \lambda}$ of $\partial A$, which we obtain by concatenating all well-orderings $<_\beta$ for $\beta < \sigma$.
    Let $a_i$ with $i < \lambda$. Then there exists a minimal $\beta < \sigma$ such that $a_i \in \nbly(A_\beta)$. 
    It follows that
    $A_\beta \subseteq  A \cup \Set{a_j \colon j < i}$ and  $ |N(a_i)|=|N(a_i) \cap (A \cup \Set{a_j \colon j < i}|$.
\end{proof}

Let $h$ be a list-colouring of an induced subgraph of $G$.
We say a vertex $a \in \dom(h)$ \emph{is happy with} $h$ if $
|\Set{v \in \dom(h) \cap N(a) \colon h(v) = h(a)}|\leq |\Set{v \in \dom(h) \cap N(a) \colon h(a) \neq h(v)}|.
$

\begin{lemma}[{\cite[Lemma 2]{shelah1990graphs}}] \label{lemma_2}
Let $G$ be a graph and suppose that $A,B^* \subseteq V(G)$ are disjoint sets of vertices such that $B^*$ is $A$--\saturated. %\ \textcolor{blue}{and $A$, $B^*$ disjoint?}. % satisfy the conditions of the \textcolor{purple}{previous} lemma. 
Suppose further that  every $v \in V$ has a list $L(v)$ of $3$ colours, and let $C = \bigcup_{v \in V}L(v)$.
Then whenever $h \colon A \cup B^* \to C$ is a list-colouring, then there is a list-colouring $g \colon \overline{A \cup B^*} \to C$ extending $h$ such that every element of $B' \cup \partial (A \cup B^*)$ where 
$$B' = \set{b \in B^*}:{|N(b) \cap \partial (A \cup B^*)| > |N(b) \cap (A \cup B^*)}$$
is happy with $g$.
\end{lemma}
\begin{proof}
    Let $b \in B'$. 
    Since $|N(b) \cap \partial(A \cup B^*)| > |N(b) \cap (A \cup B^*)|$, we have $N(b) \setminus (A \cup B^*) \neq \emptyset$.
    Since $B^*$ is $A$--\saturated, we get
    \[
        |(N(b) \setminus A) \cap B^*| = |B^*|.
    \]
    Hence, we have $|N(b) \cap (A \cup B^*)| \geq |B^*|$.
    By the definition of $B'$ it follows that
    \[
    |N(b) \cap \partial(A \cup B^*)| > |B^*|.
    \]
    Hence, by Lemma \ref{lemma_1} there is a family $\Set{F_b \colon b \in B'}$ of pairwise disjoint sets such that
    \[
    F_b \subseteq N(b) \cap \partial(A \cup B^*) \quad \text{and} \quad |F_b| =| N(b) \cap \partial(A \cup B^*)|= |N(b) \cap \overline{A \cup B^*}|  
    \]
    for every $b \in B'$. The second equation follows by the definition of $B'$.
    By Proposition \ref{prop_enumerate_boundary} there exists an enumeration $\Set{z_i \colon i < \lambda}$ of the elements of $\partial(A \cup B^*)$ such that
    \[
    |N(z_i)| = |N(z_i) \cap (A \cup B^* \cup \{z_j \colon j < i\})| %\quad (i < \lambda).
    \]
    for every $i < \lambda$.
    
    We extend $ h $ to the function $g\colon \overline{A \cup B^*} \to C$ by choosing $ g(z_i) \in L(z_i) $ inductively for every $ i < \lambda $. 
    At the $i$-th step there are two possible choices for $ g(z_i) $ that will ensure that $z_i$ is happy with $g$.
    Consequently, if $ z_i \in F_b $ for some $b \in B'$, then we may also choose $ g(z_i) $ different from $ g(b) $. The function $ g $ so constructed satisfies the requirements of the lemma.
\end{proof}

\begin{proof}[Proof of Theorem~\ref{thm_intro2} for uncountable $B$] 
    
   Suppose for a contradiction that Theorem~\ref{thm_intro2} is false and consider a counterexample with $\mu = |B|$ minimal. By the results of the previous section, $\mu$ is uncountable. Let us make the following observations:
   \begin{enumerate}[label=$(\arabic*)$]
        \item\label{clm_1} \emph{For every infinite $B' \subseteq B$ with $|B'|< |B|$ we have $\overline{A \cup B'} \neq V$.}
    \end{enumerate}

    \begin{proof}[Proof of \ref{clm_1}]
    \renewcommand{\qedsymbol}{$\Diamond$}
    Suppose  there exists an infinite $B' \subseteq B$ with $\overline{A \cup B'}=V$ and $|B'| < \mu$.
    We can assume that $x \in B'$.
    By Lemma \ref{closure_lemma}, there exists $B^*$ with $B' \subseteq B^*$ and $|B'|=|B^*|$ such that $B^*$ is $A$--\saturated\ and $B^* \cap A = B' \cap A = \emptyset$.
    By the minimality of $\mu$ there exists $h' \colon A \cup B^* \to C$ extending $h$ such that every vertex of $B^*$ is happy with $h'$.
    Next, we apply Lemma~\ref{lemma_2} to extend $h'$ to a list colouring $g \colon V= \overline{A \cup B^*} \to C$ such that every vertex of $V \setminus (A \cup B^*) \cup B''$ with
    \[
        B'' = \set{b \in B^*}:{|N(b) \cap (V \setminus (A \cup B^*))| > |N(b) \cap (A \cup B^*)|}
    \]
    is happy with $g$.
    It remains to show that all vertices of $B^*\setminus B''$ are happy with $g$.
    Let $v \in B^* \setminus B''$. 
    Then we have $|N(v) \cap (V \setminus (A \cup B^*))| \leq |N(v) \cap (A \cup B^*)|$, let us call this inequality $(*)$.
    We claim that $|N(v)|=|N(v) \cap (A \cup B^*)|$.
    If $N(v) \setminus (A \cup B^*) \neq \emptyset$, then we get $|(N(v) \setminus A) \cap B^*|=|B^*|$ by $B^*$ being $A$-saturated. So $v$ has infinite degree, and since one cannot partition an infinite set into two stricty smaller sets, the claim follows from $(*)$. Otherwise, if $N(v) \setminus (A\cup B^*)=\emptyset$, then $N(v)=N(v) \cap (A \cup B^*)$ and the claim is trivial. 
    It follows that $v$ has already $d(v)$ many opposite neighbours under $h'$, and so is automatically happy with $g$.
    \end{proof}

   \begin{enumerate}[label=$(\arabic*)$, resume]
        \item\label{clm_2} \emph{For every infinite $B' \subseteq B$ with $|B'|< |B|$ we have $| B \setminus \overline{A \cup B'}| = \mu$.}
    \end{enumerate}

       \begin{proof}[Proof of \ref{clm_2}]
    \renewcommand{\qedsymbol}{$\Diamond$}
    If $R =  B \setminus \overline{A \cup B'}$ has size $<\mu$, then the infinite set $B'' = R \cup B'$ has size $<\mu$, but the closure of $A \cup B''$ is equal to $V$, contradicting \ref{clm_1}. 
    \end{proof}

    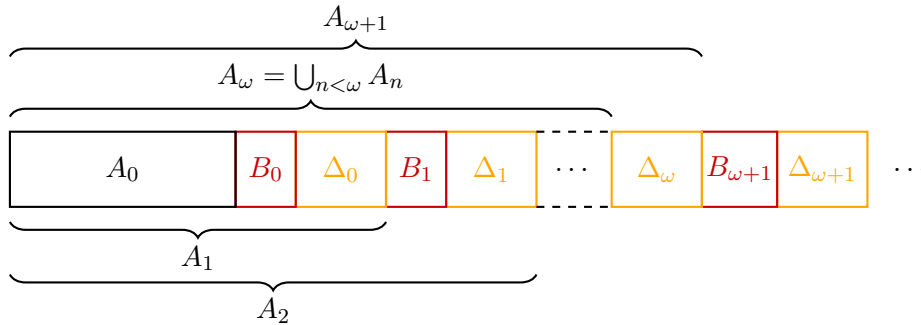
\begin{figure}[h]
        \centering
        \begin{tikzpicture}
    % Define colors
    \definecolor{orange}{RGB}{255,165,0}
    \definecolor{red}{RGB}{200,0,0}

    % Blocks (order is inversed, so that the colours overlap in the right way)
    \draw[thick, orange] (5.8,0) rectangle (7,1) node[midway, orange] {$\Delta_1$};
    \draw[thick, red] (5,0) rectangle (5.8,1) node[midway, red] {$B_1$};
    \draw[thick, orange] (3.8,0) rectangle (5,1) node[midway, orange] {$\Delta_0$};
    \draw[thick, red] (3,0) rectangle (3.8,1) node[midway, red] {$B_0$};
    \draw[thick] (0,0) rectangle (3,1) node[midway] {$A_0$};
    
    % Dots instead of Dotted Extension
    \node at (7.5,0.5) {$\dots$};
    \draw[dashed, thick] (7,0) -- (8,0);
    \draw[dashed, thick] (7,1) -- (8,1);
    
    % Next set
    \draw[thick, orange] (10.2,0) rectangle (11.4,1) node[midway, orange] {$\Delta_{\omega+1}$};
    \draw[thick, red] (9.2,0) rectangle (10.2,1) node[midway, red] {$B_{\omega+1}$};
    \draw[thick, orange] (8,0) rectangle (9.2,1) node[midway, orange] {$\Delta_\omega$};
    \node at (12,0.5) {$\dots$};
    
    % Lower Braces 
    \draw[decorate, decoration={brace, amplitude=6pt, mirror}, thick] (0,-0.2) -- (5,-0.2) node[midway, below=5pt] {$A_1$};
    \draw[decorate, decoration={brace, amplitude=6pt, mirror}, thick] (0,-0.9) -- (7,-0.9) node[midway, below=5pt] {$A_2$};
    % Upper braces
    \draw[decorate, decoration={brace, amplitude=6pt}, thick] (0,1.2) -- (8,1.2) node[midway, above=5pt] {$A_\omega = \bigcup_{n < \omega} A_n$};
    \draw[decorate, decoration={brace, amplitude=6pt}, thick] (0,2) -- (9.2,2) node[midway, above=5pt] {$A_{\omega+1}$};

\end{tikzpicture}
        \caption{The construction}
        \label{fig:enter-label}
    \end{figure}

 We now construct subsets $\set{A_\alpha}:{\alpha \leq \mu}$ and $\set{B_\alpha}:{\alpha < \mu}$ of $V$ such that

    \begin{enumerate}[label=$(\arabic*)$, resume]
        \item\label{clm_3} $A_0 =A$, $A_{\alpha+1} = \overline{A_\alpha \cup B_\alpha}$, $A_\alpha = \bigcup_{\beta < \alpha} A_\beta$ (if $\alpha$ is a limit), and $A_\mu = V$,
        \item\label{clm_4} $x \in B_0$ and $B_\alpha \subseteq  B \setminus A_\alpha$ for every $\alpha < \mu$,
        \item\label{clm_5} $B_\alpha = \emptyset$ if $\alpha$ is a limit,
        \item \label{clm_6} otherwise, if $\alpha$ is a successor, then
        \smallskip
        \begin{enumerate}[label=$(6.\arabic*)$]
            \item%[$(6.1)$]
            \label{clm_6.1} $|B_\alpha|= |\alpha| + \omega$
        \end{enumerate}    
        and $B_{\leq \alpha} \coloneqq \bigcup_{\beta \leq \alpha} B_\beta$ is $A_\alpha$--\saturated, i.e.~for every $y \in B_{\leq \alpha}$ with $N(y) \setminus (A_\alpha \cup B_\alpha) \neq \emptyset$ we have
        \smallskip
        \begin{enumerate}[label=$(6.\arabic*)$, resume]
            \item%[$(6.2)$]
            \label{clm_6.2} $|N(y) \setminus A_\alpha | > |B_{\leq \alpha}|=|B_\alpha|$ and
            \item%[$(6.3)$]
            \label{clm_6.3}  $|(N(y)\setminus A_\alpha)  \cap B_\alpha|=|B_{\leq \alpha}|=|B_\alpha|$.
        \end{enumerate}
    \end{enumerate}
    We also write $$\Delta _{\alpha}:= A_{\alpha+1} \setminus \p{A_{\alpha} \cup B_{\alpha}} =\partial (A_{\alpha} \cup B_\alpha),$$
    so $A_{\alpha+1} = A_\alpha \cup B_\alpha \cup \Delta_\alpha$. Note that we control the sizes of the $B_\alpha$, but not the sizes of the $\Delta_\alpha$.

  \begin{proof}[Construction for \ref{clm_3}--\ref{clm_6}]
    \renewcommand{\qedsymbol}{$\Diamond$}
    First, we fix a well-ordering of $B$ in order type $\mu$.
    Suppose that for some ordinal $\gamma < \mu$ we have already constructed $\set{A_\beta}:{\beta < \gamma}$ and $\set{B_\beta}:{\beta < \gamma}$. 
    If $\gamma$ is a limit, then we set $A_\gamma = \bigcup_{\beta < \gamma} A_\beta$ and $B_\gamma= \emptyset$ and are done.
    
   If $\gamma = \alpha+1$ is a successor, by \ref{clm_1} we define $A_{\alpha+1}\coloneqq \overline{A_\alpha \cup B_\alpha}$.
   Then $|B \setminus A_{\alpha+1}| = |B \setminus \overline{A \cup B_{< \gamma}}|=\mu$ by \ref{clm_6.1}, \ref{clm_2} and \ref{clm_4}.
    Therefore, we may choose an infinite subset $B' \subseteq B \setminus A_{\alpha+1}$ with $|B'|=|\gamma|+ \omega$ that contains the minimal element of $B \setminus A_{\alpha+1}$ relative to the previously fixed enumeration of $B$.
    
    Then we apply Lemma \ref{closure_lemma} to $A_{\alpha+1}$, $B_{\leq \alpha}\cup B'$ to obtain $B^*$ of size $|B_{\leq \alpha}\cup B'|=|\gamma|+\omega = |B'|$ with $(B_{\leq \alpha}\cup B') \subseteq B^*$ that is $A_{\alpha+1}$--\saturated\ and $B^* \cap A_{\alpha+1} = B_{\leq \alpha}$. We set $B_{\alpha+1}\coloneqq B^* \setminus A_{\alpha+1}$, then $B^*=B_{\leq \alpha}\sqcup B_{\alpha +1}$. 
    Hence, the set $B_{\leq \alpha +1}$ is $A_{\alpha+1}$--\saturated.
    To see that \ref{clm_6.1} is satisfied, note that $B' \subseteq B_{\alpha+1} \subseteq B^*$, so all three sets have size $|\gamma|+\omega$.
    
    Finally, since we have chosen in each construction step the minimal element of $B \setminus A_\gamma$ and since $B_\beta \subseteq A_\alpha$ for every $\beta < \alpha$, it follows that every element of $B$ will be contained in some $A_\alpha$ with $\alpha < \mu$.
    This ensures that the construction stops with $A_\mu = V$.
    \end{proof}

    For every infinite cardinal $\kappa < \mu$ let $Y_\kappa = \set{y \in B_{< \kappa^+}}:{|N(y) \cap B_{<\kappa^+}|=\kappa^+}$.

     \begin{enumerate}[label=$(\arabic*)$, resume]
        \item\label{clm_7} \emph{There are pairwise disjoint sets $I_\kappa(y) \subseteq \Set{ \alpha \colon \kappa \leq \alpha < \kappa^+}$ ($y \in Y_\kappa$) with $|I_\kappa(y)|= \kappa^+$ such that $N(y) \cap B_\alpha \neq \emptyset$ for every $\alpha \in I_\kappa(y)$.}
    \end{enumerate}

        \begin{proof}[Proof of \ref{clm_7}]
    \renewcommand{\qedsymbol}{$\Diamond$}

    Fix $y \in Y_\kappa$. Then there is a unique ordinal $\beta_y < \kappa^+$ such that $y \in B_{\beta_y}$. 
    We claim that $N(y) \cap B_\alpha \neq \emptyset$ for every successor ordinal $\alpha$ with $\max\Set{\kappa,\beta_y} \leq \alpha < \kappa^+$. 
    Suppose $N(y) \setminus (A_\alpha \cup B_\alpha)= \emptyset$ for some $\alpha < \kappa^+$. 
    Then $B_\beta \subseteq V \setminus (A_\alpha \cup B_\alpha)$ and hence $N(y) \cap B_\beta = \emptyset$ for every $\beta > \alpha$. A contradiction to $y \in Y_\kappa$.

    Since $|Y_\kappa| \leq |B_{< \kappa^+}| \leq \kappa^+$ by \ref{clm_6.1}, the assertion follows by applying 
    Lemma~\ref{lemma_1} to the collection $\cA = \set{A_y}:{y \in Y_\kappa}$ where $A_y$ consists of all successor ordinals $\alpha$ with $\max\Set{\kappa,\beta_y} \leq \alpha \leq \kappa^+$.
    \end{proof}
    
    For every non-limit $\alpha < \mu $ we choose $x_\alpha \in B_\alpha$ such that $x_0 = x$ and $x_\alpha \in N(y) \cap B_\alpha$ if $\alpha \in I_\kappa(y)$ for some $\kappa < \mu$ and $y \in Y_\kappa \cap B_{<\alpha}$ (otherwise, we choose $x_\alpha \in B_\alpha$ arbitrarily.) 
    In the last construction step, we recursively construct a transfinite sequence of functions 
    $(g_\alpha \colon A_\alpha \to C)_{\alpha < \mu}$ all extending $h$ and each other, such that for every successor ordinal $\alpha+1 < \mu$ we have 
    \begin{enumerate}[resume, label=$(\arabic*)$]
        \item \label{prop_different_colour} $g_{\alpha+1}(x_{\alpha}) \neq g_{\alpha+1}(y)$, if there are $\kappa < \mu$ and $y \in Y_\kappa \cap B_{<\alpha}$ such that  $\alpha \in I_\kappa(y)$,
        \item \label{item7} every vertex of $\Delta_\alpha \cup B'_{\leq\alpha}$ is happy with $g_{\alpha +1}$, where  
        \[
            B'_{\leq\alpha} \coloneqq \Set{ y \in B_{\leq \alpha} \colon |N(y) \cap \Delta_\alpha| > |N(y) \cap (A_\alpha \cup B_\alpha)|}. 
        \]
    \end{enumerate}
 
      \begin{proof}[Construction for \ref{prop_different_colour}--\ref{item7}]
    \renewcommand{\qedsymbol}{$\Diamond$}We begin the construction with $g_0 \coloneqq h$.
    In the limit case we define $g_\alpha = \bigcup_{\beta < \alpha} g_\beta$.
    We describe the successor step $\alpha \mapsto \alpha + 1$. So suppose that $g_\alpha \colon A_\alpha \to C$ has already been defined for some $\alpha < \mu$. We begin with the following claim.
    
    \begin{enumerate}[resume, label=$(\arabic*)$]
        \item \label{g_prime}
        \emph{$g_\alpha$ extends to $g'_\alpha$ with domain $A_\alpha \cup B_\alpha$ such that every vertex of $B_\alpha$ is happy with $g'_\alpha$.}
        \end{enumerate}
        If $\alpha$ is a limit, then $B_\alpha = \emptyset$ by \ref{clm_5} and we set $g'_\alpha \coloneqq g_\alpha$.
   Otherwise, if $\alpha$ is a successor, then $A_\alpha$ is closed by \ref{clm_3}.
    Since $|B_\alpha| < \mu$, by the induction hypothesis applied to the subgraph $G_\alpha$ induced by $A_\alpha \cup B_\alpha$ there exists a list-colouring $g'_\alpha \colon A_\alpha \cup B_\alpha \to C$ that extends $g_\alpha$ such that every vertex of $B_\alpha$ is happy. 
    Furthermore, the inductive hypothesis also allows us to arrange that $g'_\alpha (x_\alpha) \neq c_\alpha$
    where $c_\alpha \coloneqq g_\alpha(y)$ if there are $\alpha \in I_\kappa(y)$ of some cardinal $\kappa < \mu$ and $y \in Y_\kappa \cap B_{<\alpha}$;
    $c_\alpha \coloneqq c_x$ if $\alpha=0$ and otherwise we choose $c_\alpha$ as an arbitrary colour of $L(x_\alpha)$.
    \end{proof}
    Since $B_\alpha \subseteq V \setminus A_\alpha$, $B_{< \alpha} \subseteq A_\alpha$ and $B_{\leq \alpha}$ is $A_\alpha$--\saturated, it follows that $B_\alpha$ is also $A_\alpha$--\saturated. 
    Now we apply Lemma \ref{lemma_2} to $g'_\alpha \colon A_\alpha \cup B_\alpha \to C$ to obtain $g_{\alpha+1} \colon A_{\alpha+1} \to C$ that extends $g'_\alpha$ such that all required properties are satisfied. %\textcolor{blue}{CHECK}

    To complete the proof, we must show that $g \coloneqq \bigcup_{\alpha < \mu} g_\alpha$ is our desired list-colouring of $G$. By \ref{clm_1}, the domain of $g$ is all of $V(G)$, and since all $g_\alpha$ extend $h$, so does $g$. 
     First, note that we have $g(x)=g_0(x_0)\neq c_x$. 
       
    It remains to show that every vertex $v \in B$ is happy with $g$.
    If $v \not \in B_{< \mu}$, then $v \in \Delta_\alpha$ for some $\alpha < \mu$.
    By Proposition \ref{prop_enumerate_boundary} (i), we have $|N(v)|=|N(v) \cap A_{\alpha +1}|$.
    Since $v$ is happy with $g_{\alpha +1}$ by \ref{item7} and $|N(v)|=|N(v) \cap A_{\alpha +1}|$, it follows that $v$ is happy with $g$.

    Otherwise, $v \in B_{<\mu}$, so by \ref{clm_5} there is successor ordinal $\alpha < \mu$ with $v \in B_\alpha$.
    Let $\beta \leq \mu$ be minimal so that 
    $d(v) = d_{A_\beta}(v)$. 
    Note that $\beta$ exists, since $A_\mu=V$.
    Furthermore, since $A_\alpha$ is closed and $v \not \in  A_\alpha$, we have $\beta > \alpha$.
    Now, we consider the following cases.

    \begin{figure}[h]
        \centering
        \begin{tikzpicture}
    % Define colors
    \definecolor{orange}{RGB}{255,165,0}
    \definecolor{red}{RGB}{200,0,0}

    % Blocks (order is inversed, so that the colours overlap in the right way)
    \draw[thick, orange] (5.8,-1) rectangle (7,1) node[midway, orange] {};
    %\node[orange] at (6.4,0) [yshift=-0.4cm] {$\Delta_{\alpha+1}$};
    \draw[thick, red] (5,-1) rectangle (5.8,1) node[midway, red] {};
    %\node[red] at (5.4,0) [yshift=-0.4cm] {$B_{\alpha+1}$};
    \draw[thick, orange] (3.8,-1) rectangle (5,1) node[midway, orange] {};
    \node[orange] at (4.4,0) [yshift=-0.4cm] {$\Delta_{\alpha}$};
    \draw[thick, red] (3,-1) rectangle (3.8,1) node[midway, red] {};
    \node[red] at (3.4,0) [yshift=-0.4cm] {$B_{\alpha}$};
    \draw[thick] (0,-1) rectangle (3,1) node[midway] {};
    \node[black] at (1.5,0) [yshift=-0.4cm] {$A_\alpha$};

    % Vertex with neighbourhood
    \fill[blue] (3.4,0.4) circle (2pt);
    \node[blue, below] at (3.4,0.4) {$v$};
    \draw[blue, thick] (3.4,0.4) -- (10.2,0.0);
    \draw[blue, thick] (3.4,0.4) -- (10.2,0.8);
    
    % Dots instead of Dotted Extension
    \node at (7.5,0) {$\dots$};
    \draw[dashed, thick] (7,-1) -- (8,-1);
    \draw[dashed, thick] (7,1) -- (8,1);
    
    % Next set
    \draw[thick, orange] (10.2,-1) rectangle (11.4,1) node[midway, orange] {};
    \node[orange] at (10.8,0) [yshift=-0.4cm] {$\Delta_\gamma$};
    \draw[thick, red] (9.2,-1) rectangle (10.2,1) node[midway, red] {};
    \node[red] at (9.7,0) [yshift=-0.4cm] {$B_\gamma$};
    \draw[thick, orange] (8,-1) rectangle (9.2,1) node[midway, orange] {};
    \node[orange] at (8.6,0) [yshift=-0.4cm] {};
    \node at (12,0) {$\dots$};

    \draw[thick, blue] (10.2,0.4) ellipse (0.8cm and 0.4cm) node[midway, align=center] {};
    \node[blue] at (10.2,0.4)  {1.A 1.B};
    
    % Lower Braces 
    \draw[decorate, decoration={brace, amplitude=6pt, mirror}, thick] (0,-1.2) -- (3.8,-1.2) node[midway, below=5pt] {$A_{\alpha+1}$};
    % Upper braces
    \draw[decorate, decoration={brace, amplitude=6pt}, thick] (0,1.2) -- (9.2,1.2) node[midway, above=5pt] {$A_\gamma$};
    \draw[decorate, decoration={brace, amplitude=6pt}, thick] (0,2) -- (11.4,2) node[midway, above=5pt] {$A_{\beta}= A_{\gamma+1}$};

\end{tikzpicture}
  
        \caption{Case 1}
        \label{fig:enter-label1}
    \end{figure}

    \textbf{Case 1} \textit{$\beta = \gamma+1$ is a successor ordinal}.
    
    \textbf{Case 1.A} 
    $d(v) = d_{B_\gamma}(v)$.
    We show that $\gamma = \alpha$. Then we have 
$d(v) = d_{B_\alpha}(v)$
    and therefore $v$ is happy with $g'_\alpha \colon A_\alpha \cup B_\alpha \to C$ by \ref{g_prime}.
    Since $g$ extends $g'_\alpha$ , $v$ is happy with $g$.
    
   So suppose for a contradiction that $\alpha < \gamma$. Then for every successor $\zeta$ with $\alpha \leq \zeta < \gamma$
    we know from the minimality of $\beta$ that
    $v$ has at least one neighbour outside of $A_\zeta \cup B_\zeta \subset  A_\gamma$, so $|N(v) \cap B_\zeta | = |B_\zeta| = |\zeta|+\omega$ by \ref{clm_6.3}. 
But this implies \[
        |N(v) \cap A_\gamma|\geq  \sup \Set{|\zeta| + \omega \colon \alpha \leq \zeta < \gamma} = |\gamma| + \omega = |B_\gamma|\geq d(v). 
    \]
    Thus, $d(v)=d_{A_\gamma}(v)$, contradicting the minimality of $\beta$.

    \textbf{Case 1.B} 
     $d(v) > d_{B_\gamma}(v)$. By the definition of $\beta$ and $\gamma < \beta$, we have $|N(v)|> |N(v) \cap A_\gamma|$.
    Then it follows that
    \[
        d(v)= |N(v) \cap A_{\gamma+1}|=|N(v) \cap\Delta_{\gamma}| > |N(v) \cap (A_\gamma \cup B_\gamma)|.
    \]
    The first equation follows from the definition of $\beta=\gamma+1$.
    The second equation follows from $|N(v)|>|N(v) \cap A_\gamma|$ and $|N(v)|>|N(v) \cap B_\gamma|$. 
  Therefore, it follows from property \ref{item7} that $v \in B'_{\leq \gamma}$.
  With $d(v)=d_{A_{\gamma+1}}(v)$ it follows that $v$ has $d(v)$-many opposite neighbours in $g_{\gamma+1}$. 
  Hence, $v$ is happy with $g$.

  \medskip
    
    \textbf{Case 2} \textit{$\beta$ is a limit ordinal.}
    
    \textbf{Case 2.A} $d(v) \leq |\beta|$.
    First, we observe that for every successor $\zeta$ with $\alpha \leq \zeta < \beta$ 
      we know from the minimality of $\beta$ that
    $v$ has at least one neighbour outside of $A_\zeta \cup B_\zeta \subset  A_{\zeta+1} $, so $|N(v) \setminus A_\zeta| > |B_\zeta| = |\zeta| + \omega$ by  \ref{clm_6.2}.
   
    This implies $|\beta|=\beta$: Indeed, if $\zeta:= |\beta| < \beta$, then 
    $$|\beta| \geq |N(v)|\geq |N(v) \setminus A_\zeta| > |B_\zeta| = |\zeta| + \omega = |\beta| +\omega =  |\beta|,$$ a contradiction. 
    Then \ref{clm_7} implies that 
    \begin{align*}
         & |\Set{\zeta \colon \alpha \leq \zeta < \beta, \zeta \in I_\kappa(v) \text{ for some cardinal } \kappa < \beta}| \\ 
         =& \sum_{|\alpha| \leq \kappa < |\beta|} |I_\kappa(v) |  = \sum_{|\alpha| \leq \kappa < |\beta|} \kappa^+ = |\beta|.
    \end{align*}

    Hence, it follows from \ref{prop_different_colour} that $g_{\zeta+1}(v) \neq g_{\zeta+1}(x_\zeta)$ for $|\beta|$-many $\zeta < \beta$.
    Since $g_\beta$ extends every $g_\zeta$ with $\zeta < \beta$, it follows that $v$ has $d(v)$-many opposite neighbours in $g_\beta$. Therefore, $v$ is happy with $g$. 
    
    \textbf{Case 2.B} $\lambda := d(v)>|\beta|$.
    For $\zeta$ with $\alpha < \zeta < \beta$ let us write $N_\zeta = N(v) \cap \Delta_\zeta$, a set of size $<\lambda$ by the minimality of $\beta$. As $|N(v) \cap B_{<\beta}| \leq |B_{<\beta}| \leq |\beta| < d(v)$, and $|N(v) \cap A_\beta| = d(v) = \lambda$ by the definition of $\beta$, it follows that $N = \bigcup_{\alpha \leq \zeta < \beta} N_\zeta$ has size $\lambda$. 
    In particular, $\lambda$ is a singular cardinal, because the $\lambda$-sized set $N$ is a union over $|\beta| < \lambda$ many sets $N_\zeta$, all of size $<\lambda$. Furthermore, the minimality of $\beta$ implies that $|\bigcup_{\alpha \leq \zeta < \beta^*} N_\zeta| < \lambda$ for all $\beta^* < \beta$ by the minimality of $\beta$.
    
     \begin{enumerate}[resume, label=$(\arabic*)$]
        \item \label{clm_11}
        \emph{There are an increasing sequence of cardinals $(\lambda_i)_{i < cf(\lambda)}$ with $\lambda_0 > |\beta|$ and an increasing sequence of ordinals $(\beta_i)_{i < cf(\lambda)}$ such that $\lambda = \sup \lambda_i$ 
    and\[
        |N_{\beta_i}|=\lambda_i > |N(v) \cap (A_{\beta_i} \cup B_{\beta_i})|.
    \]
    }
    \end{enumerate}

            \begin{proof}[Proof of \ref{clm_11}]
    \renewcommand{\qedsymbol}{$\Diamond$}
    Fix an auxiliary  increasing sequence $(\lambda'_i)_{i < cf(\lambda)}$ of cardinals that are cofinal for $\lambda$ with $\lambda'_0 > |\beta|$. 
    Suppose that for some $j < cf(\lambda)$ we have already chosen $(\beta_i)_{i < j}$ and $(\lambda_i)_{i < j}$ correctly. Let $\beta^*_j = \sup_{i < j} \beta_i$ and  $\lambda^*_j = \sup_{i < j} \lambda_i$. Then $|N(v) \cap A_{\beta^*_j}| \leq \sup_{i < j} \lambda_i = \lambda^*_j < \lambda$, so $\beta^*_j < \beta$ by the definition of $\beta$.
    Now choose $\beta_j \geq \beta^*_j$ such that $|N_{\beta_j}|$ is larger than $\lambda'_j$ and $\lambda^*_j$, which is possible since $|\bigcup_{\alpha \leq \zeta < \beta^*_j} N_\zeta| < \lambda$ but $|\bigcup_{\alpha \leq \zeta < \beta} N_\zeta| = \lambda$.
      \end{proof}
    According to property \ref{item7}, we get $ v \in B'_{\leq \beta_i}$ for all $i < cf(\lambda)$. In particular, $v$ is happy with $g_{\beta_i+1}$, so $v$ has $\lambda_i$  neighbours of a colour different from its own under $g_{\beta_i+1}$. But from this it follows that $v$ has  $\lambda = \sup \lambda_i$  neighbours of a colour other than its own under $g_\beta$ and hence under $g$, so $v$ is happy with $g$ as desired.
\end{proof}
\bibliographystyle{plain}
\bibliography{ref}

\end{document}